\documentclass{amsart}
\usepackage{amsmath}
\usepackage{amsmath}
\usepackage{amssymb}

\input xy
\xyoption{all}

\newcommand{\ran}{\mathop{\mathrm{ran}}}

\newcommand{\free}[1]{\underset{\scriptscriptstyle #1}{\displaystyle{\ast}}\,}

\title{Certain free products of graph operator algebras}

\author{Benton L. Duncan}

\address{Department of Mathematics\\
300 Minard Hall\\
North Dakota State University\\
Fargo, ND  58105-5075\\
USA}

\email{benton.duncan@ndsu.edu}



\begin{document}

\theoremstyle{plain}
\newtheorem{thm}{Theorem}
\newtheorem{lem}{Lemma}
\newtheorem{prop}{Proposition}
\newtheorem{cor}{Corollary}

\theoremstyle{definition}
\newtheorem{dfn}{Definition}
\newtheorem*{construction}{Construction}
\newtheorem*{example}{Example}

\theoremstyle{remark}
\newtheorem*{conjecture}{Conjecture}
\newtheorem*{acknowledgement}{Acknowledgements}
\newtheorem{rmk}{Remark}

\begin{abstract}

We develop a notion of a generalized Cuntz-Krieger family of
projections and partial isometries where the range of the partial
isometries need not have trivial intersection. We associate to these
generalized Cuntz-Krieger families a directed graph, with a coloring
function on the edge set.  We call such a directed graph an
edge-colored directed graph. We then study the $C^*$-algebras and
the non-selfadjoint operator algebras associated to edge-colored
directed graphs. These algebras arise as free products of directed
graph algebras with amalgamation.  We then determine the
$C^*$-envelopes for a large class of the non-selfadjoint algebras.
Finally, we relate properties of the edge-colored directed graphs to
properties of the associated $C^*$-algebra, including simplicity and
nuclearity. Using the free product description of these algebras we
investigate the $K$-theory of these algebras.
\end{abstract}

\maketitle

The operator algebras of directed graphs are important for two
reasons. First because they give concrete examples of large classes
of operator algebras. Second, they are useful because structural
properties of the algebras can be related to simple observations
about their underlying graphs. In this paper we take this nice class
of algebras and look at their universal free products amalgamated
over specific subalgebras. Again, the reasons are two-fold, we
believe that this large class of ``concrete" examples will lead to a
better understanding of universal free products, and we also have
seen that the defined underlying discrete structure will provide
insight into structural questions concerning these algebras. The
focus of this paper, besides introducing this class of free
products, is to make a case that the second aim of this research is
tractable.

Graph algebras are a generalization of the Cuntz algebras where we
have a collection of projections, and collections of partial
isometries with domain and range satisfying natural conditions
corresponding to the collection of projections.  These relations can
be identified via directed graphs with projections corresponding to
vertices and arrows between projections corresponding to partial
isometries.  The source of the arrow corresponds to the domain
projection for the partial isometry, and the range of the arrow
corresponds to the range projection for the partial isometry.
However, in the graph algebra context we require that for any
projection corresponding to a vertex the set of edges ending at the
vertex the ranges of the associated partial isometries ``sum" to the
projection.

However, when we take free products of directed graph algebras with
amalgamation over the subalgebra generated by the projections the
restrictions on the range of the partial isometries disappear; we
can then allow partial isometries with the same range. By
amalgamating over the subalgebra corresponding to the projections we
avoid some of the technicalities in \cite{Duncan:2004} and
\cite{Duncan:2004a}, where similar free products are studied.  Some
complications do persist and this additional complexity gives rise
to a more complicated discrete structure since we need to keep track
of which ranges sum to a particular projection and which ones do not
have this property. To deal with this we add a coloring to the edge
set of the directed graph which allows us to keep track of when
partial isometries have interacting range projections.  We have
called such a graph with a coloring function on the edge set a
edge-colored directed graph.  In this paper we introduce both the
$C^*$-algebras and the non-selfadjoint operator algebras associated
to edge-colored directed graphs.

Following the discussion of edge-colored directed graphs and
discussion of the edge-colored directed graph $C^*$-algebras we
prove an analogue of gauge invariant uniqueness. We then investigate
the non-selfadjoint algebras.  We then look at structural properties
starting with simplicity, and then discussing nuclearity and
exactness. Using a free product result of Thomsen
\cite{Thomsen:2003} we discuss the $K$-theory for edge-colored
directed graph $C^*$-algebras which is built up using free products
and the $K$-theory for $C^*$-algebras of directed graphs as
described in Chapter 7 of \cite{Raeburn:2005}. We will assume that
the reader is familiar with the theory of graph $C^*$-algebras and
point to \cite{Raeburn:2005} for an excellent introduction to the
theory.

\section{Cuntz-Krieger families and edge colored directed graphs}

An edge-colored family of partial isometries on a Hilbert space
$\mathcal{H}$ is a triple $(P,E,f)$ where $P$ is a collection of
pairwise orthogonal projections in $\mathcal{B}(\mathcal{H})$, $E$
is a collection of partial isometries in $\mathcal{B}(\mathcal{H})$,
and $f$ is a function of $E$ into the natural numbers.

\begin{dfn} We say that such an edge-colored family of partial
isometries on $\mathcal{H}$ call it $(P,E,f)$, is an edge-colored
Cuntz-Krieger family on $\mathcal{H}$ if $\{ P,f^{-1}(n) \}$ is a
Cuntz-Krieger family on $ \mathcal{H}$ for each $n \in \mathbb{N}$.
\end{dfn}

Notice that any Cuntz-Krieger family will clearly be an edge-colored
Cuntz-Krieger family via the map $ f(S) = 1$ for all $ S \in E$.  On
the other hand if $S_1$ and $S_2$ are partial isometries such that $
S_1^*S_1 = S_2^*S_2 = S_1S_1^* = S_2S_2^*$ then letting $ P =
S_1S_1^*$, $E = \{ S_1, S_2 \}$, and $ f(S_i) = i$, then $(P,E,f)$
is an edge-colored Cuntz-Krieger family which is not a Cuntz-Krieger
family.  It is obvious from the preceding example that the map $f$
is not unique.  In particular, for this example any function $g: E
\rightarrow \{ m_1,m_2\}$ with $ g(S_i) = m_i$ where the $m_i$ are
distinct positive integers will yield an edge-colored Cuntz-Krieger
family $(P, E, g)$.  We would like to not have to make a distinction
between these two examples.

\begin{dfn} If $(P,E,f)$ and $(Q,F,g)$ are edge-colored Cuntz-Krieger
families on $ \mathcal{H}$ then we say $(P,E,g)$ is equivalent to
$(Q,F,f)$ if $P = Q$ and $E = F$.
\end{dfn}

This defines an equivalence relation on edge-colored Cuntz-Krieger
families.  Further we can now assume up to equivalence that $
\ran(f) = \{ 1, 2, \cdots, m \}$ for some $ \infty \geq m \geq 1$.

\begin{dfn} If $(P,E,f)$ is an edge-colored Cuntz-Krieger family on
$\mathcal{H}$ then define $c(f) = \sup \{ \ran (f) \}$.  We define
\[ c(P,E,f) := \inf \{ c(g) : (P,E,g) \mbox{ is equivalent to }
(P,E,f) \}. \] The number $ c(P,E,f)$ will be called the coloring
number of $(P,E,f)$. \end{dfn}

As with Cuntz-Krieger families we wish to associate a combinatorial
object to edge-colored Cuntz-Krieger families.  To do this we will
need the notion of an edge-colored directed graph.  We will then see
that this association is related to the free products of graphs
described in \cite[Definition 4]{Duncan:2004}.

\begin{dfn} An edge-colored directed graph is a directed graph $G =
(V,E,r,s)$ and a function $ f: E \rightarrow \mathbb{N} $. The
function $f$ will be called a {\em edge-coloring} of $G$.  We say
that two colorings $f$ and $g$ on $G$ are equivalent if for each $v$
there is a map $\theta_v: G \rightarrow G$ which is a directed graph
automorphism such that $ \theta_v|_{V(G)}$ is the identity map, $
\theta_v$ fixes $G_v$, and there is a one-to-one map $\tau:
\mathbb{N} \rightarrow \mathbb{N}$ such that $f^{-1}(i) \cap G_v =
g^{-1}(\tau(i)) \cap G_v$.
\end{dfn}

We now establish some notation.  For each vertex $v \in V(G)$ let
$G_v$ denote the directed graph with $ V = V(G)$ and $E = \{ e \in
E(G): r(e) = v \}$. Notice that if $(G,f)$ is an edge-colored
directed graph then $(G_v,f)$ is an edge-colored directed subgraph
of $(G,f)$.

While the definition for equivalence of edge-colorings is awkward we
will see that this corresponds to the equivalence of edge-colored
Cuntz-Krieger families.  It is not difficult to see that equivalence
of edge-colorings induces an equivalence relation. We will not make
it a habit of differentiating between equivalent edge-colorings.
Again we notice that an edge-coloring of $G$, call it $f$, can be
chosen up to equivalence, with $ \ran(f) = \{ 1, 2, \cdots, m \}$
for some $ \infty \geq m \geq 1$.  In this case we say that $ f $ is
an $m$-coloring of $G$, and set $m:= c(f)$. Now define $ c(G,f) =
\inf \{ c(g): (G,g) \mbox{ is equivalent to } (G,f) \}$, and call
this the coloring number of $(G,f)$.

\begin{prop} Let $(G,f)$ denote an edge-colored directed graph.  If
$m(v)$ is the number of colors incident at $v$, then $c(G,f) =
\sup_{v \in V(G)} m(v)$. \end{prop}

\begin{proof} Notice that $m(v) \leq C(G,f)$ for all $v$, by
definition, and hence if $\sup_{v \in V(G)} m(v) = \infty$ we are
done.

So assume that there exists some vertex $v_0$ such that $ m(v_0)
\geq m(w)$ for all vertices $v \in V(G)$.  Next for each vertex $v$
order $\{ f(e): e \in G_v \} := \{ a_1, a_2, \cdots, a_m \}
\subseteq \mathbb{N}$ with the usual ordering.  Then define $f_v(e)
= \{ i: f(e) = a_i \}$.  Notice that $\sum_{v \in V(G)}f_v$ will be
a coloring of $G$ which is equivalent to $f$.  Notice further that $
c(\sum f_v) \leq m(v_0)$ and the result follows. \end{proof}

Let $\{G_{\lambda} \} = \{ ( V(G_{\lambda}), E(G_{\lambda}),
r_{\lambda}, s_{\lambda}) \}$ be a collection of directed graphs
with a common subcollection $\mathcal{V} \subseteq V(G_{\lambda})$
for all $ \lambda$.  We define the free product graph
$\free{\mathcal{V}}G_{\lambda}$ to be the graph given by:
\begin{enumerate} \item $\displaystyle{V\left(\free{\mathcal{V}}G_{\lambda}\right) = \mathcal{V}
\cup \left( \bigcup_{\lambda} V(G_{\lambda} \setminus \mathcal{V})
\right)}$.

\item $\displaystyle{E\left(\free{\mathcal{V}}G_{\lambda}\right) = \bigcup_{\lambda}
E(G_{\lambda})}$.

\item If $ e \in G_{\lambda}$ then $r(e) = r_{\lambda}(e)$.

\item If $ e \in G_{\lambda}$ then $s(e) = s_{\lambda}(e)$.
\end{enumerate}

We now present a fairly straightforward connection between
edge-colored Cuntz-Krieger families and edge-colored directed
graphs.

\begin{prop} Given an edge-colored Cuntz-Krieger family on
$\mathcal{H}$, call it $(P,E,f)$, there is an edge-colored directed
graph $(G,h)$ associated to $(P,E,f)$.\end{prop}

\begin{proof} Define a directed graph $G$ by setting $ V(G) := P $
and $ E(G) := E$.  For $S \in E$ define $ s(S) := S^*S \in P$.  Next
define $r(S) := p \in P$ such that $ SS^* \leq p$.  This completes
the description of the directed graph $G$.  Notice that $(G,f)$ will
be an edge-colored directed graph. \end{proof}

Similarly, we will see that given an edge-colored directed graph
$(G,f)$ there is always an edge-colored Cuntz-Krieger family
$(P,E,f)$ with associated edge-colored directed graph $(G,f)$.  In
fact we will show that there is a universal $C^*$-algebra for each
edge-colored directed graph.

\section{$C^*$-algebras of edge-colored directed graphs}

We wish to define a universal $C^*$-algebra of an edge-colored
directed graph.  We will start by stating the universal property we
would like it to satisfy.  We will then show using free products the
existence of such a $C^*$-algebra.

Let $(G,f)$ be an edge-colored directed graph and assume that
$(P,E,f)$ is an edge-colored Cuntz-Krieger family associated to
$(G,f)$. Notice that $ \{ p \in P \} \cup \{ s \in E \} \cup \{ s^*:
s \in E \} $ generates a $C^*$-algebra, call it $C^*(P,E,f)$.

\begin{dfn} We say that a $C^*$-algebra $A$ is universal for a
edge-colored directed graph $(G,f)$ if \begin{itemize} \item $A$ is
generated by an edge-colored Cuntz-Krieger family $(P,E,f)$
associated to $(G,f)$ and \item given any edge-colored Cuntz-Krieger
family $(Q,F,g)$ associated to $(G,f)$ there is a $*$-representation
$\pi: A \rightarrow C^*(Q,F,g)$.\end{itemize}  If such a universal
algebra exists we will call it $C^*(G,f)$.
\end{dfn}

We first establish the existence of the $C^*$-algebra.

\begin{thm}\label{existence} Given an edge-colored directed graph $(G,f)$ the algebra
$C^*(G,f)$ exists and can be written as a free product \end{thm}

\begin{proof} If we let $G_i$ denote the directed graph $(V(G),
f^{-1}(i), r, s)$ then notice that $G = \cup G_i$. Let $P_i$ denote
the collection of projections in $G_i$ associated to the vertices in
$G_i$ and notice that there is a natural $*$-isomorphism between the
$P_i$'s. We will call this subalgebra $P$ and view it as sitting
inside $C^*(G_i)$ in the natural way.  We now claim that $ C^*(G,f)
= \free{P} C^*(G_i)$. We will denote the usual Cuntz-Krieger family
for $C^*(G_i)$ by $(P,E_i)$.

We define an edge-colored Cuntz-Krieger family by $(P,\cup E_i, f)$
where $f(S) = i$ where $ S \in E_i$.  Further, by construction, the
graph associated to $(P,\cup E_i,f)$ will be $(G,f)$.  We need only
verify the universal property.  This is simply a matter of applying
universal properties for the free product.

Assume that $(Q,F,g)$ is another edge-colored Cuntz-Krieger family
associated to $(G,f)$.  Then $Q =P$ and $F = \cup E_i$ and hence
there is a $*$-representation $\pi_i: C^*(G_i) \rightarrow
C^*(P,E_i)$.  Now using the free product we have a
$*$-representation $
*\, \pi_i: \free{P}C^*(G_i) \rightarrow C^*(Q,F,g)$ which is onto since $C^*(Q,F,g)$
is generated by the set $\{ P, \cup E_i \}$.  The result now
follows.
\end{proof}

Since $C^*(G,f)$ is generated by a Cuntz-Krieger $(G,f)$ family it
will satisfy the properties for $C^*(P,E,f)$ as above.

A further consequence of the free product construction is the
following result which applies to graph algebras and hence allows
graph algebras to be written as free products of simpler algebras.

\begin{thm}\label{localityofcoloring} Let $(G,f)$ be an edge-colored
directed graph with $V(G) = \{ v_{\lambda } \}_{\lambda \in \Lambda
}$.  For $ \lambda \in \Lambda$ define $ E_{\lambda} := e \in E(G)$
such that $ r(e) = v_{\lambda}$. Then $G_{\lambda}:= (V(G),
E_{\lambda}, f)$ is an edge-colored directed graph where we view,
$f, s,$ and $r$ as restrictions of $f$, the source map, and the
range map, respectively.  Further, there is a natural embedding
$C^*(G_{\lambda,f}) \rightarrow C^*(G,f)$ such that for $P$, the
subalgebra of $C^*(G,f)$ generated by the vertex projections we have
$C^*(G,f) = \free{P} \{ C^*(G_{\lambda},f) \}$.
\end{thm}

\begin{proof} Notice that $C^*(G_{\lambda},f)$ is a
$C^*$-algebra generated by the edge-colored Cuntz-Krieger family
$(P, E_{\lambda})$.  It follows that $\free{P} \{ C^*(G_{\lambda},f)
\}$ is a $C^*$-algebra generated by the edge-colored Cuntaz-Krieger
family $(P,\cup \{E_{\lambda} \})$.  It is not difficult to see that
the edge-colored directed graph associate to this new edge-colored
Cuntz-Krieger family is $(G,f)$ and hence there is, by universality,
a $*$-representation $\pi: C^*(G,f) \rightarrow \free{P} \{
C^*(G_{\lambda},f)\}$.

For the reverse arrow notice that the natural inclusion of
$C^*(G,\lambda,f) $ into $C^*(G,f)$ induces a $*$-representation of
$\free{P} \{ C^*(G_{\lambda},f) \}$ into $C^*(G,f)$ which is onto a
generating set. It is a simple matter to verify that this
representation is the inverse of the representation $\pi$.
\end{proof}

As a corollary we get the following result about $C^*(G)$ where $G$
is a directed graph.

\begin{cor} Let $G$ be a directed graph and for every
vertex $v$ let $C^*(G_{v})$ be the subalgebra of $C^*(G)$ generated
by $ G_v$, then $C^*(G) = \free{\{ P_v: v \in V(G) \} }
\{C^*(G_{v})\}$.
\end{cor}

As an example we know that $M_2(\mathbb{C})$ is $C^*(G)$ where \[ G
= \xymatrix{ {\bullet} \ar[r] & {\bullet} }. \]  Noticing that this
is just a subgraph of \[ H = \xymatrix{ {\bullet} \ar[r] & {\bullet}
\ar[r] & {\bullet} } \] and letting $D = \left\{
\begin{bmatrix} a & 0 & 0 \\ 0 & b & 0 \\ 0 & 0 & c \end{bmatrix}
\right \}$  we have \[ \begin{bmatrix} \mathbb{C} & 0 \\
0 & M_2(\mathbb{C})
\end{bmatrix} \free{D}
\begin{bmatrix} M_2(\mathbb{C}) & 0 \\ 0 & \mathbb{C}
\end{bmatrix} = M_3(\mathbb{C}) \] since $ M_3(\mathbb{C})$ is the
$C^*$-algebra of the graph $H$.

We can put this result together with the general free product
description from Theorem \ref{existence} to see that the algebras
$C^*(G,f)$ are ``locally" copies of free products of simpler
algebras.  In particular, the simpler algebras are directed graph
algebras arising from directed graphs $G$ with $|V(G)| \leq 2$ and $
|\{ r(e): e \in E(G) \}| \leq 1$.

We are now in a position using free products to translate many of
the well known theorems for $C^*$-algebras of directed graphs into
the context of edge-colored directed graphs.  We first use free
products to develop an analogue of the Gauge Invariant Uniqueness
Theorems for $1$-colorable graphs.

\begin{thm} Let $(G,f)$ be an edge-colored directed graph with
$C^*(G,f) = \free{P} C^*(G_i)$ where $G_i$ is a $1$-colored directed
graph for each $i$.  Let $A_i$ be a collection of $C^*$-algebras
with common subalgebra $D$.  If $\pi_i: C^*(G_i)\rightarrow A_i$ is
a $*$-isomorphism for each $i$ with $ \pi_i|_P \in D$ for all $i$
and $\pi_i = \pi_j|_P$ for all $ i ,j$ then the induced
representation $*\, \pi_i: C^*(G,f) \rightarrow \free{\pi_i(P)}
\pi_i(C^*(G_i))$ is an isomorphism.\end{thm}

\begin{proof} This follows by noticing that $ \pi_i^{-1}$ will
induce a $*$-representation which will be the inverse of the
$*$-representation $*\, \pi_i$. \end{proof}

To see this as an analogue of the Gauge Invariant Uniqueness Theorem
notice that one needs only verify a gauge action on a Cuntz-Krieger
family $(P,S)$ that corresponds to the gauge action on the $G_i$ and
then $A_i = C^*(P,S)$ will be isomorphic to $C^*(G_i)$ and the above
result will apply.

Before looking at properties of the $C^*$-algebras of edge-colored
directed graphs, we wish to develop the non-selfadjoint variant of
the construction.

\section{Non-selfadjoint operator algebras of edge-colored directed
graphs}

We say that a map $\pi: G \rightarrow B(\mathcal{H})$ is a
contractive representation of $(G,f)$ if the restriction of $\pi$ to
the directed graph $(V(G), f^{-1}(j), r, s )$ is a contractive
representation for each $j \in \mathbb{N}$.

Given an edge-colored directed graph $(G,f)$ we will define the
algebra $A(G,f)$ to be the norm closed operator algebra satisfying
the universal property.

(Universal Property for $A(G,f)$) There exists a contractive
representation $\iota: G \rightarrow A(G,f)$ such that for $ \pi:G
\rightarrow B(H)$ a contractive representation of $(G,f)$ there
exists a unique completely contractive representation $ \tilde{\pi}:
A(G,f) \rightarrow B(H)$ satisfying $ \tilde{\pi} \circ \iota =
\pi$.

We focus first on existence and uniqueness of this algebra.  This
follows exactly as in the $C^*$-case by using universal free
products.

\begin{thm} Given an edge-colored directed graph $(G,f)$ the algebra
$A(G,f)$ exists and is unique. \end{thm}

\begin{proof} Again letting $G_i$ be the directed graph given by
$(V(G), f^{-1}(i), r, s)$. Denoting the graph algebra $A(G_i)$ by
$A_i$ and letting $P$ denote the subalgebra of $A(G_i)$ generated by
the vertex projections we claim that $\free{P}A_i$ will satisfy the
universal property reserved for $A(G,f)$.  Uniqueness will then
follow in the usual manner for universal objects.

So let $\pi: G \rightarrow B(\mathcal{H})$ be a contractive
representation of $(G,f)$.  Then, as $G_j$ is a subgraph of $G$ we
have that $ \pi|_{G_i} := \pi_i$ is a contractive representation of
$G_i$. By the universal property for $A(G_i)$ there exists a
contractive representation $\iota_i: G_i \rightarrow A(G_i)$ and
completely contractive representation $\widetilde{\pi_i} : A(G_i)
\rightarrow B(\mathcal{H})$ satisfying $\widetilde{\pi_i} \circ
\iota_i = \pi_i$. By construction we know that $\pi_i|_{P} =
\pi_k|_{P}$ and hence there exists a completely contractive
representation $ \free{i} \widetilde{\pi_i} : \free{P}A_i
\rightarrow B(\mathcal{H})$ extending each of the maps $
\widetilde{\pi_i}$.

Next let $ \iota: G \rightarrow \free{P}(A_i)$ be given by $
\iota(v) = \iota_i(v)$ for any $v \in V(G_i)$, and $ \iota(e) =
\iota_i(e)$ for all $e \in E(G_i) \subseteq E(G)$.  Under the
coloring $f$ we have that this is a contractive representation of
$G$ since $ f^{-1}(i) = E(G_i)$ for all $i$. Checking that $
\free{i} \widetilde{\pi_i} \circ \iota = \pi$ we have that
$\free{P}A_i$ satisfies the universal property and hence we are
done.\end{proof}

As a corollary of the proof we can actually write $A(G,f)$ as a free
product of directed graph operator algebras, exactly as with
$C^*(G,f)$. Using the notation from the preceding proof, given an
edge-colored directed graph $(G,f)$ define for each $1 \leq i \leq
\infty$ $A_i = A(G_i)$ where $G_i$ is the directed graph $(V(G),
f^{-1}(i), r,s)$. Next let $P$ denote the algebra $A(V)$ where $V$
is the directed graph $(V(G), \emptyset, r,s)$.

\begin{cor} For an edge-colored directed graph $(G,f)$ we have
$A(G,f) = \free{P} A(G_i)$. \end{cor}

As was done for graph algebras \cite{Katsoulis-Kribs:2004} we state
a dilation theorem. This is similar to
\cite{Davidson-Katsoulis:2007}. It will be a simple corollary of the
construction using free products. As such it will not be entirely
satisfactory.

\begin{prop} Let $\pi: (G,f) \rightarrow B(\mathcal{H})$ be a
contractive representation, then there exists a Hilbert space
$\mathcal{K} \supseteq \mathcal{H}$, a contractive representation
$\tilde{\pi} : (G,f) \rightarrow B(\mathcal{K})$ satisfying:
\begin{enumerate} \item $\tilde{\pi}(e)$ is a partial isometry for
each edge $e \in E(G)$.
\item $\pi(e) = P_{\mathcal{H}} \tilde{\pi}(e)|_{\mathcal{H}}$ for all
edges $e \in E(G)$. \end{enumerate} \end{prop}

\begin{proof} For each $i$ let $\pi_i: G_i \rightarrow \mathcal{K_i}$
denote the dilation of $ \pi :G_i \rightarrow \mathcal{H}$.  Notice
that $\mathcal{H} \subseteq \mathcal{K}_i$ for all $i$, and for each
$e \in E(Q)$ there is $i$ such that $ \pi_i(e)$ is a partial
isometry with $\pi(e) = P_{\mathcal{H}} \pi_i(e)|_{\mathcal{H}}$.

Now letting $ \mathcal{K} = \mathcal{H} \oplus \left(\oplus_{i}
(\mathcal{K}_i \ominus \mathcal{H})\right)$ and define $
\tilde{\pi}: G \rightarrow B(\mathcal{K})$ by $ \tilde{\pi}(v) =
\pi(v) \oplus (\pi_i(v) - \pi(v))$ for $v \in V(G)$ and for $e \in
E(G_i) \subseteq E(G)$ define $ \tilde{\pi}(e) = \pi_i(e)$.  Notice
that $\tilde{\pi}$ will induce a representation as described.
\end{proof}

Notice that if $V(G) = 1$ and $f(e) = f(g)$ implies $e = g$, then
the algebra $A(G,f)$ is an example of the semicrossed products for
multivariable dynamics of \cite{Davidson-Katsoulis:2007}, whose free
product description motivated \cite{Duep}.

We now use the classical notion of $C^*$-envelopes to relate the
selfadjoint algebra and the non-selfadjoint algebra.

\section{$C^*$-envelopes}

To discuss the $C^*$-envelope of the algebra $A(G,f)$ we first
remind the reader of some results from \cite{Duep}. Recall that an
operator algebra $A$ has the {\em unique extension property} if for
every faithful representation $\pi$ of $C^*_e(A)$ the only
completely contractive extension of $\pi|_A$ to all of $C^*_e(A)$ is
$\pi$.  The main result of \cite{Duep} was the fact that the class
of algebras with the unique extension property is closed with
respect to universal free products.  The following was then a
corollary of this result.

\begin{thm} (\cite{Duep}) Let $A_i$ be a collection of
operator algebras with the unique extension property.  Further,
assume that the $A_i$ share a common $C^*$-subalgebra $D$.  Then
$C^*_e(\free{D}(A_i)) = \free{D} C^*_e(A_i)$. \end{thm}

We will now show that a certain subclass of directed graph algebras
$A(G)$ have the unique extension property and hence we will be able
to, as long as we stay in this subclass, determine the
$C^*$-envelope of certain $A(G,f)$.

\begin{prop}Let $G$ be a row finite directed graph, then $A(G)$ has
the unique extension property. \end{prop}

\begin{proof} Let $ \pi: C^*(G) \rightarrow B(\mathcal{H})$ be
a faithful $*$-representation.  Notice that $\pi|_{A(G)}$ will be
completely contractive and hence any completely positive extension
of $ \pi|_{A(G)}$, call it $\tau$,  must satisfy the property $
\tau(a^*) = \tau(a)^*$ for all $ a \in A(G)$.  Now letting $S_e$ be
the generator of $A(G)$ corresponding to the edge $e \in E(G)$ we
know that $\tau(S_e^*S_e) = \tau(P_{s(e)}) = \tau(S^e)^*\tau(S_e)$
since $ P_{s(e)}$ is in $A(G)$ and corresponds to the projection
onto the domain of $S_e$.

On the other hand, we know that $ \tau(S_eS_e^*) \geq
\tau(S_e)\tau(S_e)^* = \pi(S_e) \pi(S_e)^*$.  Since $G$ is row
finite then for any vertex $v$ we have \[ \sum_{\{e : r(e) = v \}} =
\pi(S_e)\pi(S_e^*).\] Now $ \pi(S_e) \pi(S_e)^*$ is a projection
orthogonal to $ \pi(S_f) \pi(S_f^*)$ for all $ f \neq e$ so we have
\[ \sum_{\{ e : r(e) = v \}} \tau(S_eS_e^*) \geq \sum_{\{ e : r(e) =
v \} } \pi(S_eS_e^*) : = \pi(P_{v}) \] with equality if and only if
$ \tau(S_eS_e^*) = \pi(S_eS_e^*)$.  But we know that $\pi(P_v) =
\tau(P_v)$ since $P_v \in A(G)$ and hence we have that $
\tau(S_eS_e^*) = \tau(S_e) \tau(S_e^*)$.  It follows that for all $
e \in E(G)$ the operator $S_e $ is in the $*$-subalgebra on which $
\tau$ is a $*$-representation, \cite[Theorem 2.18]{Paulsenbook}. But
$\{ S_e \} $ generates $C^*(G)$ since $G$ is row-finite and hence $
\tau$ is a $*$-representation and is equal to $\pi$.  Thus $A(G)$
has the unique extension property. \end{proof}

Recall that an edge-colored directed graph is said to be row finite
if $G_j$ is row-finite for all $ 1 \leq j \leq \infty$.

\begin{thm} If $(G,f)$ is a row finite edge-colored directed graph,
then $C^*_e(G,f) = C^*(G,f)$. \end{thm}

\begin{proof} Since the set of algebras with the unique extension
property is closed with respect to free products, it follows that
$A(G,f) = \free{P}A(G_i)$ has the unique extension property. It
follows that $C^*_e(A(G,f)) = \free{P} C^*_e(A_i) = \free{P}
C^*(A_i) = C^*(G,f)$ by combining the free product description for
$C^*(G,f)$ with \cite[Theorem 1]{Duep}.
\end{proof}

Of course this question remains open in the case of the non-row
finite edge-colored directed graphs.  In particular we have the
following example due to Muhly and Solel \cite{Muhly-Solel:1999}.

{\bf Example:} $A_{\infty}$ does not have the unique extension
property:

This follows the discussion in the last two paragraphs of
\cite{Muhly-Solel:1999}. In particular if $\{ t_i \}_{i=1}^{\infty}$
is a collection of partial isometries in $B(H)$ such that $ 1_{B(H)}
- \sum_{i=1}^{\infty} t_it_i^* > 0$, then there is a faithful
completely contractive representation of $O_{\infty}$ given by $ S_i
\mapsto t_i$ which is not a boundary representation for $A_{\infty}
\subset O_{\infty}$.  Of course we know that $C^*_e(A_{\infty}) =
\mathcal{O}_{\infty}$ and it seems natural to conjecture that
$C^*_e(A(G,f)) = C^*(G,f)$ in the not row-finite case.  We do not,
as yet, have any evidence to support this conjecture.

We now return to looking at the $C^*$-algebras of edge-colored
directed graphs with a focus on structural properties as well as
examining their $K$-theory.

\section{Properties of $C^*$-algebras of edge-colored directed graphs}

Many of these results will illustrate the differences between
$C^*$-algebras of edge-colored directed graphs and the usual
$C^*$-algebras of directed graphs. As with directed graph algebras
we can find graph theoretic conditions describing simplicity.

\begin{prop} If $C^*(G,f) = \free{P} C^*(G_i)$ is simple then
$C^*(G_i)$ is simple for all $i$ and given any pair of vertices $
v_1, v_2$ the subgraph $((\{ v_1,v_2 \}, E_{v_1,v_2}, r, s),f)$ is
$1$-colorable.
\end{prop}

\begin{proof} We will use the contrapositive.  First notice that
if any of the $C^*(G_i)$  are not simple then without loss of
generality there is $I \subseteq C^*(G_1)$ which is a nontrivial
ideal.  Then there is a nontrivial representation using the
universal property for free products of $C^*(G,f)$ onto
$(C^*(G_1)/I) \free{P} \left( \{ \free{P} C^*(G_i): i \geq 2 \}
\right)$ and hence $C^*(G,f)$ is not trivial.

Similarly assume that there exists vertices $v_1$ and $v_2$ such
that the subgraph with vertex set $ \{ v_1, v_2 \}$ and the edge set
given by $ \{ e: r(e) = v_1, s(e) = v_2 \}$ is not $1$-colored.
Given the free product description of our algebras it will suffice
to assume that this subgraph is $G$ and that $G$ is $2$-colorable.
So assume that there are $m$ edges of one color, $\{ e_1, e_2,
\cdots , e_m \}$ and $n$-edges of another color $\{ f_1, f_2,
\cdots, f_n \}$ and assume without loss of generality that $m \leq
n$. Notice that $ \{ S_{f_1}, S_{f_2}, \cdots,S_{f_{m-1}},
\sum_{i=m}^{n} S_{f_i} \}$ is a Cuntz-Krieger family with $m$
partial isometries and hence there is a $*$-representation of
$C^*(\mathcal{L}_m)$ into $C^*(\mathcal{L}_n)$.  There is then a
representation of the free product algebra $C^*(G,f)$ onto
$C^(S_{f_1}, S_{f_2} \cdots, S_{f_n})$ and notice that $ S_{f_1}-
S_{e_1}$ is in the kernel of this representation.  It follows that
$C^*(G,f)$ is not simple.\end{proof}

Unfortunately the converse is not true as one can see that the
edge-colored directed graph $G$ given by 
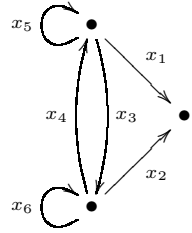
\[ \xymatrix{ {\bullet}
\ar[dr]^{x_1} \ar@/^/[dd]^{x_3} \ar@(dl,ul)^{x_5} &
\\ & {\bullet} \\ {\bullet} \ar@/^/[uu]^{x_4} \ar@(dl,ul)^{x_6} \ar[ur]_{x_2} & }
\] with coloring
\[ f(x_i) = \begin{cases} 1 & i \mbox{ is odd} \\ 2 & i \mbox{ is
even}\end{cases}.
\] is such that $C^*(G,f)$ is isomorphic to the non-simple directed
graph algebra $C^*(H)$ where $H$ is given by 
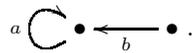
\[ \xymatrix{ {\bullet}
\ar@(dl,ul)^a  & {\bullet} \ar[l]^b }.\]  This isomorphism is
induced by the following map on generators: \begin{align*}
\tau(S_{x_1}) & \mapsto S_a S_a^* S_b^* \\
\tau(S_{x_2}) & \mapsto S_b S_b^* S_b^* \\
\tau(S_{x_3}) & \mapsto S_a S_b^* \\
\tau(S_{x_4}) & \mapsto S_a S_a^* \\
\tau(S_{x_5}) & \mapsto S_a S_a S_a^* \\
\tau(S_{x_6}) & \mapsto S_a S_b S_b^* \end{align*}

It would be interesting to find a characterization of those
edge-colored directed graphs $(G,f)$ which give rise to simple
$C^*$-algebras.  One issue that may help in the resolution of this
is in identifying when a given edge-colored directed graph algebra
is isomorphic to a directed graph algebra via a map similar to that
described in our counterexample to the converse.

We can also take known results and use them to discuss nuclearity
and exactness for these algebras.  In the $1$-colored case we know
that the algebras are all nuclear. Adding colors to the graph
complicates issues however.  Again we are limited by not having a
clear indication of when a directed graph algebra is isomorphic to a
edge-colored directed graph algebra. On the other hand, we do have
some graph theoretic characterizations of when an edge-colored
directed graph gives rise to a non-exact $C^*$-algebra.

\begin{prop} $C^*(G,f)$ is not exact if $(G,f)$ contains a subgraph
with one vertex which is not 1-colorable. \end{prop}

\begin{proof} We begin by showing that the edge-colored directed graph
$C^*$-algebra for the graph with a single vertex, and two edges of
different color is not exact.  Of course this algebra is isomorphic
to $C(\mathbb{T})\free{\mathbb{C}} C(\mathbb{T})$.  It is certainly
well known that the universal free product $C[0,1] \free{\mathbb{C}}
C[0,1]$ is not exact, for a proof see \cite{Blecher-Duncan:2007}. We
will embed this free product into $C(\mathbb{T}) \free{\mathbb{C}}
C(\mathbb{T})$ as a subalgebra. The result will then follow as
exactness is preserved by taking subalgebras.

So for $ f \in C[0,1]$ define $\widehat{f}(e^{i\theta}) =
\begin{cases} f(\frac{\theta}{\pi}) & 0 \leq \theta \leq \pi \\
f(\frac{2\pi - \theta}{\pi}) & \pi \leq \theta \leq 2 \pi
\end{cases}.$  It is a simple matter to see that this induces a
unital $*$-representation of $C[0,1]$ into a subalgebra of
$C(\mathbb{T})$.   Applying a result of Pedersen
\cite{Pedersen:1999} we get that the natural free product map
induced by $ f \mapsto \widehat{f}$ injects $C[0,1]\free{\mathbb{C}}
C[0,1]$ into $C(\mathbb{T}) \free{\mathbb{C}} C(\mathbb{T})$. The
result now follows for this algebra.

We next notice that $ C(\mathbb{T})\free{\mathbb{C}} C(\mathbb{T})$
has a natural embedding into $ \mathcal{O}_n \free{\mathbb{C}}
\mathcal{O}_m$ for all $ 1 \leq m,n \leq \infty$.  Hence if $G$ is
any graph with a single vertex and $f$ is not a $1$-coloring for $G$
then $(G,f)$ is not exact.

The final case is shown by embedding an appropriate copy of
$\mathcal{O}_n* \mathcal{O}_m$ into $C^*(G,f)$ by sending
$\mathcal{O}_n$ into the subgraph as in the hypotheses.  This will
be a faithful representation by applying \cite[Proposition
2.4]{Armstrong-Dykema-Exel-Li:2003}. Hence we have an injection of a
non-exact $C^*$-algebra into $C^*(G,f)$ and thus $C^*(G,f)$ can not
be exact. \end{proof}

Another ``forbidden'' sub-graph for uniqueness is the 3-colored
directed graph given by two vertices and three edges all of which
have range the first vertex and source the second vertex.

\begin{prop} $C^*(G,f)$ is not exact if in $(G,f)$ there exists
two vertices $v_1$, $v_2$ such that the subgraph $ ((\{ v_1,v_2 \},
\{ e: s(e) = v_1, r(e) = v_2 \},r,s),f)$ is not $2$-colorable.
\end{prop}

\begin{proof} We begin by looking at the case of the edge-colored directed graph
with two distinct vertices $v_1$ and $v_2$ and three differently
colored edges $e,f,g$ each with range $v_1$ and source $v_2$, call
this graph $(G_3,f_3)$. Notice that $ \{ S_e, S_f^*, S_g\}$ together
with the vertex projections is an edge-colored Cuntz-Krieger family
with associated graph given by reversing the arrow on $f$ in the
original directed graph.  Now this graph falls into the category of
the previous proposition and hence $C^*( S_e,S_g,S_f^*)$ is not
exact. But this algebra is the same as $C^*(S_e,S_f,S_g)$ and hence
in this case we have that the algebra is not exact.

Now if there exists vertices $v_1$ and $v_2$ and three collections
of edges given by $ \mathcal{E}:=\{ e_1, e_2, \cdots, e_n \},
\mathcal{F}:=\{ f_1, f_2, \cdots, f_m \}, $ and $ \mathcal{G} :=
\{g_1, g_2, \cdots, g_p\}$ such that each of these edges has source
$v_1$ and range $v_2$ and further $ (\{ P_{v_1}, P_{v_2}\}, \{
S_{e_i} \})$, $(\{ P_{v_1}, P_{v_2}\}, \{ S_{f_i } \})$, and $ (\{
P_{v_1}, P_{v_2}\}, \{ S_{g_i } \})$ are Cuntz-Krieger families.

In the case that the collections $ \mathcal{E}, \mathcal{F},$ and $
\mathcal{G}$ are finite we define $ \pi: C^*(G_3,f_3) \rightarrow
C^*(G,f)$ by sending $ S_e \mapsto \sum S_{e_i}, S_f \mapsto \sum
S_{f_i}, S_g \mapsto S_{g_i}$.  This will induce a
$*$-representation of $C^*(G_3,f_3)$ which is an embedding into
$C^*(G,f)$ and hence $C^*(G,f)$ can not be exact.

For the case where any of the $ \mathcal{E}, \mathcal{F}$, or $
\mathcal{G}$ are infinite we just need to choose finite subsets $
\mathcal{E}' \subseteq \mathcal{E}, \mathcal{F} ' \subseteq
\mathcal{F},$ and $ \mathcal{G}' \subseteq \mathcal{G}$ such that
\[ P:= \left(\sum_{e \in \mathcal{E}'} S_eS_e^*\right) \left(\sum_{f
\in \mathcal{F}'} S_fS_f^* \right) \left(\sum_{g \in \mathcal{G}'}
S_gS_g^* \right) \neq 0. \]  We then notice that $ (\{ P_{v_1},
P\}\{ S_eP: e \in \mathcal{E}' \}), (\{ P_{v_1}, P\}\{ S_fP: f \in
\mathcal{F}' \}),$ and $ (\{ P_{v_1}, P\}\{ S_gP: g \in \mathcal{G}'
\})$ are Cuntz-Krieger families and as in the finite case there is
an embedding of $C^*(G_3,f_3)$ into $C^*(\{ P_{v_1},P, \{ S_{e}: e
\in \mathcal{E}' \}, \{ S_{f}: f \in \mathcal{F}' \}, \{ S_{g}: g
\in \mathcal{G}' \}) \subseteq C^*(G,f)$ and again we have that
$C^*(G,f)$ is not exact. \end{proof}

There are examples of $n$-colorable graphs which are nuclear.  For
example, let $G_n$ be the graph with $n+1$ vertices $\{ v_0, v_1,
v_2, \cdots, v_n \}$ and $n$ edges $\{ e_1, e_2, \cdots, e_n \}$
such that $ r(e_i) = v_0$ for all $i$ and $s(e_i) = v_i$.  Next
define an $n$-coloring by $f(e_i) = i$. Letting $G_n^t$ denote the
graph with vertex set $\{ w_0, w_1, \cdots, w_n\}$ and edge set $
\{f_1, f_2, \cdots f_n \}$ satisfying $r(f_i) = w_i$ and $ s(f_i) =
w_0$ for all $i$, we define a $1$-coloring on $G_{n}^t$ by setting
$g(e_i) = 1$. Notice that since $G_n^t$ has a $1$-coloring we know
that $ C^*(G_n^t,g)$ is nuclear. But it is a simple exercise to see
that the map $S_e \mapsto S_f^*$ induces an isomorphism between
$C^*(G_n,f)$ and $C^*(G_n^t,g)$ and hence $C^*(G,f)$ is nuclear.

Of course one would like a complete classification of
nuclearity/exactness using only combinatorial properties.  It seems
that the presence of distinct differently colored cycles in the
undirected subgraph plays an important role in nuclearity, but a
good conjecture is not yet evident.  At the same time we notice that
in our two results we have that not only are the algebras
non-nuclear, but they are not exact.  It would be interesting to
know if there are non-nuclear edge-colored directed graph algebras,
which are exact.

Although we know that nuclearity is not preserved by free products
we do know that the subalgebra $\{P_v: v\in V(G)\}$ is commutative.
Now commutative $C^*$-algebras are nuclear so that the six term
exact sequence for $K$-groups for the amalgamated free product of
$C^*$-algebras described in \cite{Thomsen:2003} applies. In
particular, we will write $C^*(G,f)$ as $\free{P} C^*(G_i)$ where
$C^*(G_i)$ will be the $C^*$-algebra corresponding to a $1$-subgraph
of $G$. For a discussion of the $K$-groups of algebras of the form
$C^*(G_i)$ we refer the reader to \cite{Drinen-Tomforde:2002}.  The
following is just a reconstruction of Theorem 6.4 from
\cite{Thomsen:2003} in our context.

\begin{prop}  For a $2$-colored directed graph $(G,f)$ the following
six term exact sequence of $K$-groups is valid: \[ \xymatrix{ K_0(P)
\ar[r] & \oplus K_0(C^*(G_i)) \ar[r] & K_0(C^*(G,f)) \ar[d] \\
K_1(C^*(G,f)) \ar[u] & \oplus K_1(C^*(G_i)) \ar[l] & K_1(P) \ar[l]},
\] where $C^*(G,f) = C^*(G_1) \free{P} C^*(G_2)$.
\end{prop}

If $(G,f)$ is $1$-colored the previous is vacuous. Now since $K_1(P)
= 0$ and $ K_0(P) = \oplus_{v \in V(G)} \mathbb{Z}$, then for
$2$-colored directed graphs the Proposition reduces to
\[ 0 \rightarrow K_1(C^*(G_i)) \rightarrow K_1(C^*(G,f)) \rightarrow
\oplus_{v \in V(G)} \mathbb{Z} \rightarrow \oplus K_0(C^*(G_i))
\rightarrow K_0(C^*(G,f)) \rightarrow 0.\]

As with most exact sequences of $K$-groups the interesting thing,
and the complications, come from the connecting maps.  As an example
we can look at the case of the two colored directed graph $(G,f)$
with $1$ vertex, $2$ red edges, and $m$ blue edges. Next, let
$(G',f)$ be the two colored directed graph with $1$ vertex, $2$ red
edges, and $2$ blue edges, and $(G'',f)$ be the two colored directed
graph with $1$ vertex, $2$ red edges and $2$ blue edges.  Of course
we know that $C^*(G,f) = \mathcal{O}_2 \free{\mathbb{C}}
\mathcal{O}_n$, $C^*(G',f) = \mathcal{O}_2 \free{\mathbb{C}}
\mathcal{O}_2$, and $C^*(G'',f) = \mathcal{O}_2 \free{\mathbb{C}}
C(\mathbb{T})$.  A ``folklore" result shown to me by Bruce Blackadar
shows us that $C^*(G,f)$ is isomorphic to $C^*(G'',f)$. We proceed
inductively by showing that $ \mathcal{O}_2 \free{\mathbb{C}}
\mathcal{O}_{m-1}$ is isomorphic to $ \mathcal{O}_2
\free{\mathbb{C}} \mathcal{O}_m$.  Let $ (s_1, s_2, t_1, t_2, \cdots
, t_{m-1})$ be the generators for $\mathcal{O}_2 \free{\mathbb{C}}
\mathcal{O}_{m-1}$ and $ (s_1, s_2, t_1', t_2', \cdots, t_m')$ be
the generators for $ \mathcal{O}_2 \free{\mathbb{C}}
\mathcal{O}_{m}$.  Now define $ \phi: \mathcal{O}_2
\free{\mathbb{C}} \mathcal{O}_{m} \rightarrow \mathcal{O}_2
\free{\mathbb{C}} \mathcal{O}_{m-1}$ by $ \phi(s_i) = s_i$,
$\phi(t_i') = t_i$ for $ 1 \leq i \leq m-2$, $\phi(t_{m-1}') =
t_{m-1} s_1$, and $ \phi(t_{m}') = t_{m-1} s_2$.  Similarly, let $
\psi: \mathcal{O}_2 \free{\mathbb{C}} \mathcal{O}_{m-1} \rightarrow
\mathcal{O}_2 \free{\mathbb{C}} \mathcal{O}_{m}$ by $ \psi(s_i) =
s_i$, $ \psi(t_i) = t_i'$ for $ 1 \leq i \leq m-2$ and $
\psi(t_{m-1}) = t_{m-1}'s_1^* + t_m's_2^*$.  It is easy to see that
$ \phi$ and $ \psi$ are inverses of each other and hence induce an
isomorphism.  Notice that repeated applications of this construction
will yield the desired isomorphisms.

On the other hand examining the $K$-group exact sequences we get \[
0 \rightarrow K_1(C^*(G,f)) \rightarrow \mathbb{Z} \rightarrow 0
\oplus \mathbb{Z}/(m-1) \rightarrow K_0(C^*(G,f)) \rightarrow 0,
\]
\[ 0 \rightarrow K_1(C^*(G',f)) \rightarrow
\mathbb{Z} \rightarrow 0 \oplus 0 \rightarrow K_0(C^*(G',f))
\rightarrow 0,
\] and \[ 0 \rightarrow K_1(C^*(G'',f))
\rightarrow \mathbb{Z} \rightarrow 0 \oplus \mathbb{Z} \rightarrow
K_0(C^*(G'',f)) \rightarrow 0.\]  The second short exact sequence
allows us to compute that $K_1(C^*(G',f)) = \mathbb{Z}$ and $
K_0(C^*(G',f)) = 0$, and hence the same is true for $C^*(G,f)$. An
alternate version of this computation is due to Blackadar.

To deal with graphs with more than $2$-colors we have a technically
more complicated exact sequence but the presence of the $0$ group as
$K_1(P)$ allows us to prove the following.  We will first set some
notation.  If $(G,f)$ is an $m$-edge-colored directed graph with
$(G,f) = (G_1) \free{V} (G_2) \free{V} \cdots \free{V} (G_m)$ where
$G_i$ a $1$-colored subgraph of $(G,f)$. Let $(G_k,f_k)$ denote the
free product $(G_1) \free{V} (G_2) \free{V} \cdots \free{V} (G_k)$
as a $k$-colored sub-edge-colored directed graph of $(G,f)$.

\begin{prop} Let $(G,f)$ be an $m$-edge-colored directed graph so
that $C^*((G,f)) = C^*(G_1) \free{P} C^*(G_2) \free{P} \cdots
\free{P} C^*(G_m)$ with $G_i$ a $1$-colored subgraph of $(G,f)$.  We
have the following $6m$-cyclic exact sequence \[ \xymatrix@R=.5cm{ 0
\ar[r] & K_1(C^*(G_1) \oplus C^*(G_2)) \ar[r] & K_1( C^*(G_2,f_2))
\ar[d] \\ K_0( C^*(G_2,f_2))\ar[d] & K_0(C^*(G_1) \oplus C^*(G_2)) \ar[l] & K_0(P) \ar[l] \\
0 \ar[r] & K_1(C^*(G_2,f_2) \oplus C^*(G_3)) \ar[r] &
K_1(C^*(G_3,f_3)) \ar[d]\\  K_0( C^*(G_3,f_3)) \ar[d] &
K_0(C^*(G_2,f_2) \oplus C^*(G_2)) \ar[l]  & K_0(P) \ar[l]  \\ \vdots
\ar[d] & \vdots & \vdots \\ 0 \ar[r] & K_1(C^*(G_{m-1},f_{m-1})
\oplus C^*(G_m)) \ar[r] & K_1(C^*(G_m,f_m)) \ar[d]
\\ K_0( C^*(G_m,f_m)) \ar@/^6pc/ [uuuuuu] & K_0(C^*(G_{m-1},f_{m-1}) \oplus
C^*(G_m)) \ar[l] & K_0(P) \ar[l]}.\] \end{prop}

\begin{proof} This is just the concatenation of the short exact
sequence of Thomsen \cite{Thomsen:2003}, valid since $K_1(P) = 0$.
\end{proof}

Similarly we can apply \cite{Brownexact} to see that $Ext(C^*(G,f))$
is a group and compute it using short exact sequences as with the
$K$-groups.  We leave the details to the reader.

\bibliographystyle{plain}

\end{document}